\documentclass[12pt]{article}
\usepackage{amsmath}
\usepackage{latexsym}
\usepackage{amssymb}
\usepackage{a4wide}
\newcommand{\qed}{\Box}
\newtheorem{thm}{Theorem}[section]

\newtheorem{Defn}[thm]{Definition}
\newtheorem{Remark}[thm]{Remark}
\newtheorem{Example}[thm]{Example}
\newtheorem{prop}[thm]{Proposition}

\newenvironment{defn}{\begin{Defn}\rm}{\end{Defn}}
\newenvironment{rem}{\begin{Remark}\rm}{\end{Remark}}
\newenvironment{example}{\begin{Example}\rm}{\end{Example}}
\newenvironment{proof}{{\noindent\bf Proof.}}%
                  {\nopagebreak\hspace*{\fill}$\qed$\medskip\par}   
\newcommand{\Punkt}{\nopagebreak\hspace*{\fill}$\qed$}
\newcommand{\wb}{\overline}
\newcommand{\tensor}{\otimes}
\newcommand{\mto}{\mapsto}

\newcommand{\isom}{\cong}
\newcommand{\N}{{\mathbb N}}
\newcommand{\R}{{\mathbb R}}
\newcommand{\bF}{{\mathbb F}}
\newcommand{\bL}{{\mathbb L}}
\newcommand{\bS}{{\mathbb S}}
\newcommand{\bO}{{\mathbb O}}
\newcommand{\Q}{{\mathbb Q}}
\newcommand{\Z}{{\mathbb Z}}
\newcommand{\C}{{\mathbb C}}
\newcommand{\K}{{\mathbb K}}
\newcommand{\cU}{{\mathcal U}}
\newcommand{\cS}{{\mathcal S}}
\newcommand{\cL}{{\mathcal L}}
\newcommand{\cA}{{\mathcal A}}

\newcommand{\sub}{\subseteq}
\newcommand{\GL}{\mbox{{\rm GL}}}
\newcommand{\Lip}{\mbox{{\rm Lip}}}
\newcommand{\id}{\mbox{{\rm id}}}
\newcommand{\sbull}{{\scriptscriptstyle \bullet}}
\newcommand{\Diff}{\mbox{{\rm Diff}}}
\newcommand{\ev}{\mbox{{\rm ev}}}

\newcommand{\dl}{{\displaystyle\lim_{\longrightarrow}}}
\begin{document}
\begin{center}
{\bf\Large Aspects of \boldmath{$p$}-Adic Non-Linear
Functional Analysis}\\[4mm]
{\bf Helge Gl\"{o}ckner}\vspace{1.4mm}
\end{center}
\begin{abstract}
\noindent
The article provides an introduction to
infinite-dimensional
differential calculus over topological fields
and surveys some of its applications,
notably in the areas of infinite-dimensional
Lie groups and dynamical systems.\vspace{4mm}
\renewcommand{\thefootnote}{\fnsymbol{footnote}}
\footnotetext{\hspace*{-6.1mm}{\em Classification\/}:
26E30, 
46S10; 
22E65, 
26E15, 
26E20, 
37D10, 
46T20, 
47J07, 
58C15, 
58C20\\ 
{\em Kewords\/}:
Ultrametric calculus, non-archimedian analysis,
non-linear functional analysis, fixed point theorem,
inverse function theorem, implicit function theorem,
Lie group, invariant manifold, stable manifold,
pseudo-stable manifold, non-linear
mapping, dependence on parameters}
\renewcommand{\thefootnote}{\arabic{footnote}}
\end{abstract}
\begin{center}
{\bf\Large Introduction}
\end{center}
We describe aspects
of non-linear functional
analysis over topological fields,
considered as
the study of non-linear
mappings between topological
vector spaces,
their fixed points
and differentiability properties.
Our first aim is
to give an introduction to
the differential calculus of smooth and $C^k$-maps
between topological vector spaces over topological
fields developed in~\cite{Ber}.
This approach
generalizes Schikhof's
single variable calculus over complete
ultrametric fields (as in~\cite{Sch})
to multi- and infinite-dimensional situations.
In the case of mappings between real
locally convex spaces,
it is equivalent
to the usual locally convex calculus
(Keller's $C^k_c$-theory,
as in \cite{Kel} or \cite{Mil}).
Our second aim is to survey various recent
applications of differential calculus over topological fields,
and the specific techniques
underlying them.
In particular, we discuss the following topics:
\begin{itemize}
\item
The existence of fixed points
and their $C^k$-dependence on parameters (as in \cite{newimp});
\item
An implicit function theorem
for $C^k$-maps from arbitrary topological vector spaces
over valued fields to Banach spaces (established in~\cite{IMP}
and \cite{newimp});
\item
The construction of the main types of
infinite-dimensional Lie groups
(linear Lie groups, mapping groups, diffeomorphism groups,
direct limit groups)
over topological fields (carried out in~\cite{ZOO},
also \cite{DIR} and \cite{FUN});
\item
The construction of
invariant manifolds around hyperbolic fixed points
for dynamical systems modelled
on Banach spaces over valued fields
(as in \cite{INV} and \cite{IV2}),
by an adaptation of Irwin's method
(developed in \cite{Irw} and \cite{Ir2},
also \cite{LaW} and \cite{Wel});
\item
Some special tools of calculus
used in the preceding
constructions (results ensuring
differentiability properties
of non-linear mappings between function spaces~\cite{ZOO}
and spaces of sequences~\cite{INV};
exponential laws for function spaces~\cite{ZOO}).\vspace{2mm}
\end{itemize}
While analytic mappings between
Banach spaces over complete valued fields
and the corresponding analytic
Banach-Lie groups
are classical objects of study (see
\cite{BoF}, \cite{Bou};
also \cite{Ser} for the finite-dimensional
case), we are interested just as well in
mappings between general
topological vector spaces,
which have hardly been investigated so far.
This is essential for infinite-dimensional
Lie theory, since many interesting groups
cannot be modelled on Banach spaces.
It is also essential
to work with smooth (and $C^k$-) maps,
rather than analytic ones.
In fact, already in the real case,
typical examples of infinite-dimensional
Lie groups
(like diffeomorphism groups)
are smooth Lie groups,
but cannot be given an
analytic Lie group structure.
Even worse,
over each local field of positive characteristic,
finite-dimensional smooth Lie groups
exist which do not admit an analytic Lie group
structure compatible with their
topological group structure~\cite{NOA}.
For this reason,
we shall use smooth maps
as the basis of
Lie theory.\footnote{Another problem is that a convincing
notion of analytic map between (non-normable)
locally convex spaces over ultrametric
fields does not seem to exist yet,
in contrast to the real and complex cases.}\\[2.5mm]
We mention
that many
important techniques usually applied
in finite-dimensional non-archimedian
analysis do not possess infinite-dimensional
counterparts:
Neither the
techniques of algebraic geometry, which
help to analyze the most prominent examples of
finite-dimensional $p$-adic Lie groups
(linear algebraic groups);
nor the techniques
of rigid analytic geometry
(see \cite{BGR}, \cite{FaP}),
the non-archimedian analogue
of complex geometry and
function theory.
However, the ideas of
\emph{real} differential calculus
turn out to be extremely robust:
surprisingly large parts
of ordinary differential calculus carry over to
mappings between open subsets of topological vector spaces
over
topological fields,
once they are reformulated
in an appropriate way.
Schikhof's textbook \cite{Sch} bears
witness of this phenomenon
in the case of
functions of a single variable
(see also \cite{AaS}, \cite{Bar}, \cite{DeS}, \cite{DSm},
\cite{Kat}, \cite{KaK}, \cite{KKP} and \cite{Wsm},
part of which include functions of several variables).
The present article is based on the infinite-dimensional
extension of Schikhof's ultrametric calculus
developed in~\cite{Ber} (cf.\ \cite{Lu8} and \cite{Lud}
for a different, earlier approach to $p$-adic
infinite-dimensional differential calculus,
which will be discussed in Remark~\ref{einordng}\,(e) here).
We hope to illustrate that the techniques of real
differential calculus remain powerful
also when dealing with mappings between
topological vector spaces over topological fields.\\[2.5mm]
Our illustrating examples are
taken from Lie theory and dynamical systems.
In the theory of non-archimedian
dynamical systems,
a wealth of results is available
for one-dimensional polynomial
or analytic systems,
inspired by classical complex dynamics
(see, e.g., \cite{AaV},
\cite{Ben}, \cite{KhN}, \cite{Lub}).
Motivated by Lie-theoretic applications,
we here describe complementary results dealing
with higher-dimensional (and
infinite-dimensional) dynamical systems,
which have hardly been investigated
before. Also the study of
infinite-dimensional Lie groups
over topological fields
is of recent origin.
It started with discussions of 
irreducible representations
and invariant measures for certain infinite-dimensional
$p$-adic groups (see \cite{Lu8}, \cite{Lud}
and the references therein).
In the meantime, all major constructions
of real infinite-dimensional Lie groups
have been adapted to
more general topological fields (\cite{ZOO}, \cite{FUN}).\\[2.5mm]
Outside the realms of non-linear functional analysis,
Lie theory and
dynamical systems\linebreak
discussed in this article,
the differential calculus over topological fields
(and suitable\linebreak
commutative topological rings)
has been used to give an essentially
algebraic approach to differential geometry,
which is entirely based on differentiation
and does not involve integrals,
nor flows~\cite{Be2}.
Jordan theoretic applications
have been explored in~\cite{BaN}.
\section{Differential calculus over topological fields}
The basic idea of differential calculus
over topological fields is to call a map $C^1$
if one can pass from directional difference
quotients to directional derivatives
continuously.
This idea, and its implications,
will be described in more detail now.
Throughout this section, $\K$
denotes a topological field (i.e.,
a field, equipped with a non-discrete
Hausdorff topology which turns the field
operations into continuous mappings).
Topological $\K$-vector spaces
are defined as in the real case, and are assumed
Hausdorff. As the default, $E$, $F$ and $E_1, E_2,\ldots$
denote topological $\K$-vector spaces,
and $U\sub E$ an open subset.
To define continuous differentiability,
let $E$ and $F$ be topological $\K$-vector
spaces, and $f\colon  U\to F$ be a map
on an open subset $U\sub E$.
Then the directional difference quotient
\[
f^{]1[}(x,y,t)\; :=\;
\frac{f(x+ty)-f(x)}{t}
\]
makes sense for all $(x,y,t)$ in the subset
$U^{]1[} := \{(x,y,t)\in U\times E\times \K^\times
\colon \, x+ty\in U\}$
of $E\times E\times \K$.
To incorporate directional derivatives,
we must enlarge this set by allowing also
the value $t=0$.
Hence, we consider now
\[
U^{[1]}\;:=\; \{(x,y,t)\in U\times E\times \K
\, \colon \, x+ty\in U\}\,.
\]
Thus $U^{[1]}=U^{]1[}\cup (U\times E\times \{0\})$, as a disjoint union.
\begin{defn}\label{def1}
$f\colon U\to F$ is called \emph{continuously
differentiable} (or $C^1$)
if $f$ is continuous
and there exists a continuous map $f^{[1]}\colon U^{[1]}\to F$
which extends
$f^{]1[}\colon U^{]1[}\to F$.
\end{defn}
Thus, we assume the existence of a continuous map
$f^{[1]}\colon U^{[1]}\to F$ such that
\[
f^{[1]}(x,y,t)\; =\;
\frac{f(x+ty)-f(x)}{t}\quad\mbox{for all $(x,y,t)\in U^{[1]}$
such that $t\not=0$.}
\]
Given a $C^1$-map $f\colon  U\to F$ as before,
its directional derivative
\[
df(x,y)\,:=\, D_yf(x)\,:=\,
\lim_{0\not=t\to 0} {\textstyle
\frac{1}{t}(f(x+ty)-f(x))}
\,=\,
\lim_{0\not=t\to 0} f^{[1]}(x,y,t)
\,=\,
f^{[1]}(x,y,0)
\]
at $x\in U$ in the direction~$y\in E$
exists, by continuity of~$f^{[1]}$.
The map $df\colon U\times E\to F$ is continuous, being a partial map
of~$f^{[1]}$, and it can be shown that
the ``differential'' $f'(x):=
df(x,\sbull)\colon  E\to F$
of $f$ at~$x\in U$ is a continuous $\K$-linear map
\cite[Proposition~2.2]{Ber}.
Since $f^{[1]}(x,y,0)$ is a limit of difference quotients,
$f^{[1]}$ is uniquely determined
by~$f$.
\begin{example}
To prove that a given map is $C^1$,
usually one first writes down
directional difference quotients and then tries
to guess a continuous extension.
To illustrate this strategy,
let us show that every continuous linear map
$\lambda\colon E\to F$ between topological $\K$-vector spaces
is continuously differentiable.
For $x,y\in E$ and $t\in \K^\times$, we have
\begin{equation}\label{RHS}
\frac{\lambda(x+ty)-\lambda(x)}{t}\, =\, \lambda(y)\,,
\end{equation}
exploiting the linearity of~$\lambda$.
We now note that the right hand side
of (\ref{RHS}) makes sense just as
well for $t=0$, and defines a continuous
function
\begin{equation}\label{linf1}
\lambda^{[1]}\colon E\times E\times\K\to F\,,\qquad
\lambda^{[1]}(x,y,t):=\lambda(y)\,.
\end{equation}
Thus $\lambda$ is $C^1$, with $\lambda^{[1]}$
as just described (which is again a continuous linear
map). In particular, $\lambda'(x)=\lambda$
for each $x\in E$.
\end{example}
The same idea can be used to prove the Chain
Rule for $C^1$-maps.
\begin{prop}
If $f$ and $g$ are composable
$C^1$-maps, then also $f\circ g$ is $C^1$,
and
\[
(f\circ g)^{[1]}(x,y,t)\;=\; f^{[1]}\big(g(x), \, g^{[1]}(x,y,t), \, t\big)\,.
\]
In particular,
$d(f\circ g)(x,y)\;=\; df\big(g(x),\, dg(x,y)\big)$.
\end{prop}
\begin{proof}
For $t\not=0$, we have
\[
\frac{f(g(x+ty))-f(g(x))}{t}
= \frac{f\big(g(x)+t\, \frac{g(x+ty)-g(x)}{t}\big)-f(g(x))}{t}
= f^{[1]}\big(g(x),\,
g^{[1]}(x,y,t), \, t\big).
\]
Since the final expression makes
sense just as well for $t=0$
and defines a continuous function,
we see that $f\circ g$ is $C^1$,
with $(f\circ g)^{[1]}$ as asserted.
\end{proof}
Since $U^{[1]}$ is an open subset
of the topological $\K$-vector space $E\times E\times\K$,
it is possible to define $C^k$-maps by recursion.
\begin{defn}
Given an integer $k\geq 2$,
a map $f\colon E\supseteq U\to F$ is called
$C^k$ if it is $C^1$ and $f^{[1]}\colon U^{[1]}\to F$
is $C^{k-1}$.
As usual, $f$ is called $C^\infty$
(or \emph{smooth}) if it is $C^k$
for each $k\in \N$.
\end{defn}
The present recursive definition of $C^k$-maps
is particularly well-adapted to inductive
arguments.
In \cite{Lu8} and \cite{Lud},
a different approach to higher
differentiability was used.
\begin{example}
We have seen above that $\lambda^{[1]}$
is again continuous linear if $\lambda\colon E\to F$
is a continuous linear map.
Hence a straightforward induction shows
that $\lambda$ is $C^k$ for each $k\in \N$
and thus smooth.
A similar argument shows that
every continuous $n$-linear map
$\beta\colon E_1\times\cdots\times E_n\to F$
is smooth (use \cite[\S3.3]{Ber}
and induction).
\end{example}
\begin{rem}\label{rem1}
$C^k$-maps have many of the
properties familiar from the real
case.
\begin{itemize}
\item[(a)]
Compositions of composable $C^k$-maps
are $C^k$ (see \cite[Proposition~4.5]{Ber}).
\item[(b)]
Higher differentials can be defined
and have the usual properties:
If $f\colon E\subseteq U\to F$
is $C^k$, then the iterated directional derivative
\[
d^jf(x,y_1,\ldots, y_j):= (D_{y_j}\cdots D_{y_1}f)(x)
\]
exists for all $j\in \N$ such that $j\leq k$
and all $x\in U$, $y_1,\ldots,y_j\in E$.
The mapping\linebreak
$d^jf\colon U\times E^j\to F$
so defined is continuous,
and $d^jf(x,\sbull)\colon E^j\to F$
is a continuous, symmetric, $j$-linear
map, for each $x\in U$ (see \cite[Chapter~4]{Ber}).
\item[(c)]
Finite order Taylor expansions are available
in the following form: If $f\colon E\supseteq U\to F$
is $C^k$, with $k$ finite,
then there are
functions $a_j\colon U\times E\to F$
for $j\in \{0,1,\ldots, k\}$
and a continuous map
$R_k\colon U^{[1]}\to F$
with $R_k(x,y,0)=0$ for all $(x,y)\in U\times E$,
such that
\[
f(x+ty)\;=\;\sum_{j=0}^k t^ja_j(x,y)+t^kR_k(x,y,t)\quad
\mbox{for all $(x,y,t)\in U^{[1]}$.}
\]
%
The functions $a_0,\ldots, a_k$ and $R_k$
are uniquely determined.
Furthermore, $a_j$
is $C^{k-j}$, homogeneous of degree~$j$
in the second argument,
and $j!a_j(x,y)=d^jf(x,y,\ldots, y)$.
If $\K$ has characteristic~$0$,
we can divide by $j!$ and recover Taylor's formula
in its familiar form
with $a_j(x,y)=\frac{1}{j!}d^jf(x,y,\ldots, y)$
(see \cite[Chapter~5]{Ber} for details).
\item[(d)]
Despite its appearance, being $C^k$
is a local property: If $(U_i)_{i\in I}$
is an open cover of $U$ and $f|_{U_i}$
is $C^k$ for each $i\in I$, then
$f\colon U\to F$ is $C^k$
(see \cite[Lemma~4.9]{Ber}).
\end{itemize}
\end{rem}
The preceding facts make it possible
to define $C^k$-manifolds in the usual way.
\begin{defn}
A \emph{$C^k$-manifold}
modelled on a topological
$\K$-vector space~$E$
is a Hausdorff topological
space~$M$, equipped with a set
$\cA$ of homeomorphisms
$\phi\colon U\to V$
from open subsets of~$M$ onto open subsets
of~$E$, such that the domains $U$ cover~$M$
and the transition maps $\phi\circ \psi^{-1}$
are $C^k$ on their domain,
for all $\phi,\psi\in \cA$.
\end{defn}
Now $C^k$-maps between $C^k$-manifolds
can be defined as usual
(such that the $C^k$-property can be tested
in charts). Also
the direct product of two manifolds
can be defined as usual.
\begin{defn}
A \emph{Lie group} over a topological field~$\K$
is a group~$G$, equipped with a smooth
manifold structure modelled on a topological
$\K$-vector space which turns
the group inversion $G\to G$
and the group multiplication $G\times G\to G$
into smooth mappings.
\end{defn}
$C^k$-Lie groups for finite $k$ are
defined analogously.
\begin{rem}\label{einordng}
The $C^k$-maps considered
here are related to traditional concepts
as follows.
\begin{itemize}
\item[(a)]
First of all,
a map $f\colon \R^n\to\R^m$ is $C^k$
in our sense if and only if it is $C^k$ in the usual
sense. In fact, if $f$ is an ordinary
$C^1$-map, then
\[
f^{[1]}(x,y,t)\; :=\, \int_0^1  df(x+ty,y)\, dt
\]
defines a map $f^{[1]}\colon \R^n\times\R^n\times \R \to\R^m$
which is continuous as a parameter-dependent integral,
and produces the desired directional difference
quotient for $t\not=0$,
by the Mean Value Theorem.
Conversely, existence of $f^{[1]}$ easily entails
that $f$ is $C^1$ in the usual sense.
For higher $k$, one argues by a simple
induction.
\item[(b)]
More generally, a map $f\colon E\supseteq U\to F$
to a real locally convex space~$F$
is $C^k$ in our sense if and only if it is a ``Keller
$C^k_c$-map,'' i.e. $f$ is continuous,
the iterated directional
derivatives $d^jf(x,y_1,\ldots, y_j)$
(as in Remark~\ref{rem1}\,(b))
exist for all $j\in \N$ such that $j\leq k$,
and define continuous maps $d^jf\colon U\times E^j\to F$
(see \cite[Proposition~7.4]{Ber}).
\item[(c)]
Every $k$ times continuously
Fr\'{e}chet differentiable map ($FC^k$-map) between
real Banach spaces is $C^k$,
and conversely every $C^{k+1}$-map
is $FC^k$ (cf.\ (b) and \cite{Kel}).
\item[(d)]
A map $E\supseteq U\to F$ to a complex
locally convex space is $C^\infty$ if and only
if it is complex analytic, i.e., $f$ is continuous
and given locally around each point by a pointwise
convergent series of continuous homogeneous
polynomials \cite[Proposition~7.7]{Ber}.
\item[(e)]
Finally, a map $\K\supseteq U\to \K$
on an open subset of a complete ultrametric
field~$\K$ is $C^k$ in our sense
if and only if it is a $C^k$-map
in the sense of \cite[Definition~29.1]{Sch},
as usually considered in non-archimedian
analysis (see \cite[Proposition~6.9]{Ber}).
More generally, our $C^k$-maps
of several variables coincide with the
usual ones,
as in \cite[\S84]{Sch} and \cite{DeS}
(see \cite{COM}).
By contrast,
being $C(k)$
for $k\geq 2$
in the sense considered in \cite[\S2.3]{Lu8}
and \cite{Lud}
is a properly weaker property in general
than being $C^k$ in our
sense,
already for functions of one
variable (which therefore
do not coincide with Schikhof's):
see \cite{COM}.
\linebreak
Examples show that such $C(k)$-maps need not
admit $k$-th order Taylor expansions.
\end{itemize}
\end{rem}
Of course, we are mainly interested
in differential calculus over
a valued field $(\K,|.|)$;
thus $\K$ is a field
and $|.|\colon \K\to [ 0,\infty[$
an ``absolute value''
with properties analogous to the modulus
of real or complex numbers.
We shall always assume that $|.|$
is non-trivial in the sense
that the metric
$d(x,y):=|x-y|$ defines a non-discrete
(field) topology on~$\K$.
A valued field $(\K,|.|)$ is
called an \emph{ultrametric field}
if its absolute value $|.|$
satisfies the ultrametric inequality,
viz.\ $|x+y|\leq\max\{|x|,|y|\}$
for all $x,y\in \K$.
Besides $\R$ and $\C$,
the most prominent examples
of valued fields
are local fields (i.e., totally disconnected,
locally compact topological fields).
Any such is known to admit an ultrametric absolute
value defining its topology.
Up to isomorphism,
every local field either
is a finite extension of the field
$\Q_p$ of $p$-adic numbers,
or a field $\bF(\!(X)\!)$
of formal Laurent series
over a finite field $\bF$
(see \cite{Wei}).
If $(\K,|.|)$ is a valued field,
then we can speak of norms
and seminorms on $\K$-vector spaces,
as in the real and complex cases.
A normed space $(E,\|.\|)$
over $\K$ which is complete in the metric
determined by $\|.\|$ is called
a \emph{Banach space}.
If $(\K,|.|)$ is an ultrametric field here
and $\|.\|$ satisfies
the ultrametric inequality,
then $(E,\|.\|)$ is called an \emph{ultrametric}
Banach space.
See also \cite{BTV},
\cite{Mon} and \cite{Roo}.
\begin{rem}
We mention
that a certain range of strengthened differentiability
properties is available
for mappings between topological
vector spaces over a
valued field
$(\K,|.|)$.
In this case, also $k$ times strictly differentiable
mappings ($SC^k$-maps) and $k$ times Lipschitz
differentiable mappings ($LC^k$-maps) can be
defined \cite[\S2 and \S3]{newimp}.
Then
\[
C^{k+1}\; \Rightarrow\;
LC^k\;
\Rightarrow\;
SC^k\;
\Rightarrow
\;
C^k
\]
(see \cite[Remark~3.17]{newimp}).
The strengthened differentiability
properties
are useful
for some refined results.
Strict differentiability at a point for mappings on
normed spaces is a classical concept (see \cite[1.2.2]{BoF}).
Mappings
between real Banach spaces are strictly
differentiable
if and only if they are once
continuously Fr\'{e}chet differentiable~\cite[2.3.3]{BoF}.
If $\K$ is locally compact
and $E$ a finite-dimensional $\K$-vector
space,
then a map $f\colon U\to F$ on an open subset
$U\sub E$ is $C^k$ if and only if
it is $SC^k$ (cf.\ \cite[Lemma~3.11]{newimp}).
\end{rem}
\section{Fixed points, inverse and implicit functions}
In this section, we present
an implicit function theorem
for $C^k$-maps from arbitrary
topological vector spaces over valued fields
to Banach spaces.
We also describe some of the tools used in
its proof, notably a theorem
ensuring the $C^k$-dependence
of fixed points on parameters,
which may be of independent interest.
The results are taken from~\cite{newimp}
(for slightly less general versions,
see \cite{IMP}).
Similar implicit function
theorems in the real
and complex cases,
in different settings of analysis,
were also obtained
in \cite{Hil}, \cite{Hi2}
(for Keller's $C^k_\Pi$-theory)
and \cite{Tei} (in the Convenient
Setting of Analysis, \cite{FaK}, \cite{KaM}).\\[2.5mm]
While classical implicit function theorems
are restricted to mappings between Banach spaces,
the following result
(a special case of \cite[Theorem~5.2]{newimp})
only requires that the range space
be a Banach space (for the case $k=1$, see
Remark~\ref{nowk1}).
\begin{thm}[Generalized Implicit Function Theorem]\label{geninv}
Let $\K$ be a valued field,
$E$ be a topological $\K$-vector space,
$F$ be a Banach space over~$\K$,
and $f\colon U\times V\to F$ be a
$C^k$-map, where $U\sub E$
and $V\sub F$ are open subsets and $k\in \N\cup\{\infty\}$,
$k\geq 2$.
Assume that $f(x_0,y_0)=0$ for some $(x_0,y_0)\in U\times V$
and $f_{x_0}'(y_0)\in \GL(F)=\cL(F)^\times$,
where $f_x:=f(x,\sbull)
\colon V\to F$ for $x\in U$.
Then there exist open neighbourhoods
$U_0\sub U$ of~$x_0$ and
$V_0\sub V$ of $y_0$ such that
$\{(x,y)\in U_0\times V_0\colon f(x,y)= 0 \}$
is the graph of a $C^k$-map
$\lambda \colon U_0\to V_0$.
\end{thm}
\textbf{Strategy of the proof.}
The idea is to deduce the implicit
function theorem from an
``Inverse Function Theorem with Parameters.''
Since
$f_{x_0}'(y_0) \in \GL(F)$,
where $\GL(F)$ is open in $\cL(F)$
and the mapping $U \to \cL(F)$, $x \mto f'_x(y_0)$
is continuous, we see that
$f_x'(y_0)\in \GL(F)$ for $x$ close to $x_0$.
Then $f_x^{-1}$ exists locally
around $f_x(y_0)$,
by a suitable Inverse Function Theorem
(e.g., \cite[1.5.1]{BoF} or
Proposition~\ref{newton} below).
Now the essential point is that
the map
\begin{equation}\label{needdomain}
\psi\, \colon \, (x,z) \, \mto \, (f_x)^{-1}(z)
\end{equation}
actually makes sense on a whole
neighbourhood $U_0\times W $ of $(x_0,0)$ in
$U\times F$,
and is $C^k$ there (further explanations
will be given below).
Once we have this, the rest is easy:
The map $\lambda \colon U_0\to F$,
$\lambda (x):= (f_x)^{-1}(0)=\psi(x,0)$
is $C^k$,
and $f(x,\lambda(x))=0$.\vspace{2.5mm}\Punkt

\noindent
To see that the map $\psi$ in (\ref{needdomain})
is defined on a whole neighbourhood,
one exploits that suitable
versions of the inverse
function theorem provide
\emph{quantitative} information on the size
of the domain of the local inverse.
For example,
one can use the following Lipschitz
version of the inverse function theorem,
\cite[Theorem~5.3]{newimp}
(adapted from the real case in~\cite{Wel}).
If $E$ is a Banach
space over $\K$ and $f\colon E\supseteq X\to E$ a map,
we set
$\Lip(f) := \sup\{\|f(y)-f(x)\|\cdot \|y-x\|^{-1}\colon
x\not=y\in X\}$
and call $f$ \emph{Lipschitz} if $\Lip(f)<\infty$.
\begin{prop}[Lipschitz Inverse Function Theorem]\label{newton}
\,Let $(E,\|.\|)$ be a Banach\linebreak
space over a valued field
$(\K,|.|)$.
Let $x\in E$, $r>0$
and $f\colon  B_r(x)\to E$ be a map
on the ball $B_r(x):=\{y\in E\colon \|y-x\|<r\}$
of the form
\[
f\; =\; A+\tilde{f}\, ,
\]
where $A\in \GL(E)$
and $\tilde{f}\colon B_r(x)\to E$
is a Lipschitz map
with $\Lip(\tilde{f})<\|A^{-1}\|^{-1}$.
Then
$f$ has open image and is a homeomorphism
onto its image. Furthermore,
the inverse map $f^{-1}\colon f(B_r(x))\to B_r(x)$
is Lipschitz with
$\Lip(f^{-1})\leq (\|A^{-1}\|^{-1}-\Lip(\tilde{f}))^{-1}$,
and
\[
B_{a r}(f(x)) \; \sub \; f(B_r(x)) \; \sub \; B_{b r}(f(x))
\]
with
$a :=\|A^{-1}\|^{-1}-\Lip(\tilde{f})$ and
$b :=\|A\|+\Lip(\tilde{f})$.
\end{prop}
The proof of Proposition~\ref{newton} is based on a simple Newton-type
iteration.
Accordingly, $\psi(x,z)$ in (\ref{needdomain})
is the fixed point of a contraction $g_{x,z}$
of a (subset of a) Banach space.
The following general result
ensures that
the fixed point $\psi(x,z)$ of $g_{x,z}$
is a $C^k$-function of $(x,z)$.
It applies to so-called ``uniform'' families
of contractions.
\begin{defn}\label{defunicon}
Let $F$ be a Banach space over
a valued field~$\K$,
and $U\sub F$ be a subset.
A family $(f_p)_{p\in P}$ of mappings
$f_p\colon U\to F$
is called a
\emph{uniform family of contractions}
if there exists $\theta\in [0,1[$
such that $\|f_p(x)-f_p(y)\|\leq \theta\|x-y\|$
for all $x,y\in U$ and $p\in P$.
\end{defn}
We now formulate the technical backbone of the generalized implicit
function theorem
(a special case of \cite[Theorem~4.7]{newimp}).
It was stimulated by the discussion
of fixed points in real Banach spaces
and their dependence
on parameters in \emph{Banach}
spaces in~\cite{Ir2}.
\begin{thm}[on the Parameter-Dependence of Fixed Points]
Let $\K$ be a valued field,
$E$ be a topological $\K$-vector space,
$F$ be a Banach space over~$\K$
and
$f\colon P\times U\to F$
be a $C^k$-map,
where $P\sub E$
and $U\sub F$ are open
and $k\in \N_0\cup\{\infty\}$.
Assume that the maps $f_p:=f(p,\sbull)\colon
U\to F$ define a
is a uniform family of contractions
$(f_p)_{p\in P}\,$.
Then
\[
Q\, :=\,
\{p\in P\colon \mbox{$f_p$ has a fixed point $x_p$}\}
\]
is open in~$P$,
and $\phi\colon Q\to F$, $p\mto x_p$
is a $C^k$-map.
\end{thm}
Exploiting the $C^k$-dependence
of fixed points on parameters,
one also obtains the
following analogue of
the classical Inverse Function Theorem
(see \cite[Theorem~5.1]{newimp}):
\begin{thm}[Inverse Function Theorem]\label{basicinv}
Let $E$ be a Banach space
over a valued field~$\K$
and $f\colon U\to E$ be a $C^k$-map
on an open subset $U\sub E$,
where $k\in \N \cup\{\infty\}$
and $k\geq 2$.
If $f'(x)\in \GL(E)$ for some $x\in U$,
then there exists an open neighbourhood
$V\sub U$ of~$x$ such that $f(V)$ is open
in~$E$ and $f|_V\colon V\to f(V)$
is a $C^k$-diffeomorphism.
\end{thm}
In the ultrametric case,
closer inspection of the mapping $\psi$
in (\ref{needdomain}) leads
to the following useful result
(which is a special case
of \cite[Theorem~5.17]{newimp}).
\begin{thm}[Ultrametric Inverse Function
Theorem with Parameters]\label{invparultra}
\,Let $F$\linebreak
be an ultrametric Banach space
over an ultrametric field $(\K,|.|)$
and $E$ be a topological $\K$-vector
space.
Let $P_0 \sub E$ and
$U\sub F$ be open, and $f\colon  P_0 \times U \to F$ be a $C^k$-map,
where $k\in \N\cup\{\infty\}$ and $k\geq 2$.
Let $(p_0,x_0)\in P_0\times U$ be given such that
$A:=f_{p_0}'(x_0)
\in \GL(F)$,
where $f_p:=f(p,\sbull)\colon  U \to F$ for $p\in P_0$.
Then there exists an open neighbourhood $P\sub P_0$ of~$p_0$
and $r>0$ such that $B:=B_r(x_0)\sub U$ and the following holds:
\begin{itemize}
\item[\rm (a)]
$f_p(B)=f(p_0,x_0)+A.B_r(0)=:V$, for each $p\in P$,
and $\phi_p\colon  B\to V$, $\phi_p(y):=f(p,y)$ is a
$C^k$-diffeomorphism;
\item[\rm (b)]
$f_p(B_s(y))=f_p(y)+A.B_s(0)$
for all $p\in P$, $y\in B$ and $s\in \,]0,r]$;
\item[\rm (c)]
The map $\psi\colon  P\times V\to B$, $\psi(p,v):=\phi_p^{-1}(v)$
is $C^k$;
\item[\rm (d)]
$\xi \colon  P\times B\to P\times V$,
$\xi(p,y):=(p,f(p,y))$
is a $C^k$-diffeomorphism,
with inverse given by $\xi^{-1}(p,v)=(p,\psi(p,v))$.
\end{itemize}
\end{thm}
A similar result is available for non-ultrametric
Banach spaces \cite[Theorem~5.13]{newimp}.
\begin{rem}\label{nowk1}
Theorems~\ref{geninv} and~\ref{invparultra}
remain valid for $k=1$ if $F$ is finite-dimensional,
and Theorem~\ref{basicinv}
remains valid for finite-dimensional~$E$
(see \cite{newimp}).
In the infinite-dimensional case,
$C^1$-maps need not be approximated
good enough by their linearization
around a given point to guarantee the hypotheses of
Proposition~\ref{newton}.
This problem disappears
for $SC^k$-maps
or $LC^k$-maps.
All of the results
presented so far in this section
have analogues for $SC^k$-maps
and $LC^k$-maps,
valid for all $k\in \N\cup\{\infty\}$ (see \cite{newimp}).
\end{rem}
{\bf Applications.}
We mention some applications
of the preceding results
in non-archimedian analysis.
For applications in the real and complex
cases, see \cite{IMP} and \cite{GaN}.
\begin{itemize}
\item
As described in more
detail in Section~\ref{secinv},
the implicit function theorem
can be used to construct stable
manifolds around hyperbolic fixed points.
\item
In \cite{ZOO}, the ultrametric
inverse function theorem with parameters
in a Fr\'{e}chet space
is used to prove that the inversion map
$\Diff(M)\to \Diff(M)$, $\gamma\mto \gamma^{-1}$
of the diffeomorphism group of a paracompact,
finite-dimensional smooth manifold over a
local field is smooth
(see also Section~\ref{seclie} below).
\item
Varying an idea from~\cite{CaH},
$C^{k+n}$-solutions to (systems of)
$p$-adic differential
equations of the form $y^{(k)}=f(x,y,y',\ldots, y^{(k-1)})$
for $f$ an $LC^n$-map and $k\in \N$, $n\in \N_0$
can be constructed
using our inverse function theorems (with and without
parameters),
which depend on initial conditions and parameters
in a controlled way
(work in progress;
cf.\ \cite[\S\,65]{Sch} for $C^1$-solutions to scalar-valued
first order equations).
This is a first step towards the study
of equations which are not accessible
by the existing highly-developed techniques
for linear or analytic equations
(cf.\ \cite{Dwo}, \cite{PaS}).
\item
The ultrametric inverse function theorem
with parameters is also used in the construction
of a $C^k$-compatible analytic Lie group structure
on each finite-dimensional $p$-adic $C^k$-Lie group
in \cite{ANA} (for $k\in \N\cup\{\infty\}$).
\end{itemize}
\section{Invariant manifolds around fixed points}\label{secinv}
We briefly discuss invariant manifolds
over valued fields and an application.
\begin{defn}\label{defnhypo}
Let $E$ be a Banach space over a valued
field $(\K,|.|)$, and $a\in \;]0,1]$.
A $($bicontinuous$)$ linear automorphism
$\alpha\in \GL(E)$ is called \emph{$a$-hyperbolic}
if $E=E_1\oplus E_2$ for certain $\alpha$-invariant
closed vector subspaces~$E_1$,~$E_2$;
\[
\|x_1+x_2\|=\max\{\|x_1\|,\|x_2\|\}\quad\mbox{for all
$x_1\in E_1$ and $x_2\in E_2$}
\]
holds for a norm $\|.\|$
on~$E$ equivalent to the original one;
and $\alpha=\alpha_1\oplus \alpha_2$ with
$\|\alpha_1\|<a$ and $\|\alpha_2^{-1}\|^{-1}> a$.
The $1$-hyperbolic automorphisms
are simply called \emph{hyperbolic}.
\end{defn}
\begin{rem}
If $E$ is an ultrametric Banach space,
then the norm $\|.\|$
described in the definition
can be chosen ultrametric
as well.
\end{rem}
\begin{rem}
If $\K$ is a local field and $\dim_\K(E)<\infty$,
then $\alpha$ is $a$-hyperbolic if and only
if $a\not=|\lambda|$ for each eigenvalue $\lambda\in \wb{\K}$
of $\alpha\tensor_\K \id_{\wb{\K}}\in \GL(E\tensor_\K \wb{\K})$,
where $\wb{\K}$ is an algebraic closure of~$\K$
(\cite{INV}; cf.\ also \cite[proof of Lemma~3.3]{SCA}).
\end{rem}
We begin with the global Stable Manifold
Theorem (see~\cite{INV}).
\begin{thm}
Let $M$ be a $C^\infty$-manifold
modelled on a Banach space
over a valued field,
$f\colon M\to M$
be a $C^\infty$-diffeomorphism,
and $p\in M$ be a hyperbolic
fixed point of~$f$,
i.e.\ $f(p)=p$ and the
differential $T_pf\colon T_pM\to T_pM$
is a hyperbolic automorphism of
the tangent space $T_pM$.
Then
\[
W_s\; :=\; \{x\in M\, \colon \;\mbox{$f^n(x)\to p\,$
as $\,n\to\infty$}\}
\]
is an immersed $C^\infty$-submanifold of~$M$.
\end{thm}
\begin{rem}
Analogous conclusions are available for
manifolds and diffeomorphisms
which are $C^k$ (with $k\geq 2$),
resp., $SC^k$ and
$LC^k$ (with $k\in \N)$.
The theorem also holds for $\K$-analytic diffeomorphism
of $\K$-analytic manifolds,
in the sense of~\cite{BoF}
(see \cite{INV}).
\end{rem}
Using standard arguments,
the global stable manifold theorem
follows from corresponding
local results. Therefore, we now briefly
discuss the construction
of local $a$-stable manifolds,
for $a\in \,]0,1]$.
For simplicity, we only consider smooth maps,
and concentrate
on the case where
$E$ is an ultrametric Banach space over
an ultrametric field $(\K,|.|)$.\\[2.5mm]
Let $\alpha\in \GL(E)$ be $a$-hyperbolic,
say $E=E_1\oplus E_2$
with $E_1$ and $E_2$ as in Definition~\ref{defnhypo}.
Let $r>0$
and $f\colon U\to E$ be a smooth map
on the ball $U:=B_r^E(0)=U_1\times U_2$, where $U_j:=B^{E_j}_r(0)$,
such that $f(0)=0$ and $f'(0)=\alpha$.
The \emph{$a$-stable set of $f$} is defined as
\[
W_{s,a}\!:=\! \{z\in U \colon \! \mbox{$f^n(z)$ is defined
for all $n\in \N_0$,
$a^{-n}\|f^n(z)\|<r$ and $f^n(z)=o(a^n)$}\}.
\]
Then $f(W_{s,a})\sub W_{s,a}$.
If $z\in W_{s,a}$, then the orbit
$\omega:=(f^n(z))_{n\in \N_0}$ is an element
of the Banach sequence space
\[
\cS_a(E)\; :=\; \{z=(z_n)_{n\in \N_0}\in E^{\N_0}\, \colon\;
z_n=o(a^n)\}
\]
with norm $\|z\|_a:=\max\{a^{-n}\|z_n\|\colon n\in \N_0\}$.
After shrinking~$r$,
we may assume that
$\Lip(f-\alpha)< \min\{1,\|\alpha_2^{-1}\|^{-1}\}$.
Then the following holds:
\begin{thm}
$W_{s,a}$ is the graph of a $C^\infty$-map $\phi\colon
U_1\to U_2$ with $\phi(0)=0$
and $\phi'(0)=0$.
Thus $W_{s,a}$ is a $C^\infty$-submanifold
of~$E$, which is tangent to the $a$-stable subspace~$E_1$
at~$0$.
\end{thm}
The proof given in \cite{INV}
follows Irwin's method (see \cite{Irw} and \cite{Wel}).
The idea is to construct not $\phi(x)$ itself,
but the orbit $\omega$ of $(x,\phi(x))$,
which is an element of the ball
\[
\cU\; :=\; \{z\in \cS_a(E)\colon \|z\|_a<r\}\, .
\]
The orbit $\omega$
turns out to solve an equation
$\Phi(x,\omega)=0$,
for a suitable smooth map
\[
\Phi\, \colon \; U_1\times \cU\to \cS_a(E)\,,
\]
to which the implicit function theorem
can be applied.
This idea, used by Irwin in the real case,
works just as well in the present situation.
Essentially, one simply uses
Theorem~\ref{geninv}
instead of the implicit function theorem
for real Banach spaces;
the desired differentiability properties of~$\Phi$
are ensured by the next proposition (from~\cite{INV}).
\begin{prop}
Let $E$ and $F$ be Banach spaces over a valued
field,
$f\colon B_r^E(0)\to F$ be\linebreak
a mapping such that
$f(0)=0$, $a\in \;]0,1]$
and $\cU:=\{w\in \cS_a(E)\colon \|w\|_a<r\}$.
If $f$ is $C^k$, $SC^k$,
$LC^k$ with $k\in \N\cup\{\infty\}$,
resp., $\K$-analytic,
then also the map
\[
\cS_a(f) \, \colon  \; \cU\to \cS_a(F)\, ,\qquad
(x_n)_{n\in \N_0}\mto (f(x_n))_{n\in \N_0}
\]
is $C^k$, $SC^k$, $LC^k$, resp., $\K$-analytic.
\end{prop}
We mention that also Irwin's
construction of pseudo-stable manifolds
(see \cite{Ir2} and \cite{LaW})
can be adapted
to general valued fields~\cite{IV2}.
Combining these constructions,
all main types of invariant manifolds
(stable and unstable manifolds,
center-stable and center-unstable
manifolds, as well as center manifolds)
become available also over valued fields.
Using these invariant manifolds,
the following result (from~\cite{SPO})
can be obtained, which provided
the original stimulus
for the author's studies compiled
in this section.
\begin{thm}
Given $k\in \N\cup\{\infty,\omega\}$,
let $G$ be a finite-dimensional
$C^k$-Lie group
over a local field $\K$
and $\alpha\colon G\to G$ be a
$C^k$-automorphism whose contraction group
\[
U_\alpha:=\{x\in G\, \colon \, \mbox{$\alpha^n(x)\to 1$
as $n\to\infty$}\}
\]
is closed in~$G$.
Then $U_\alpha\,$, $U_{\alpha^{-1}}$
and $M_\alpha:=\{x\in G\colon \mbox{$\alpha^\Z(x)$ is relatively
compact}\}$
are closed Lie subgroups,
and the product map
$U_\alpha\times M_\alpha\times
U_{\alpha^{-1}}\to G$,
$(x,y,z)\mto xyz$
is a $C^k$-diffeomorphism onto
an open, $\alpha$-stable identity neighbourhood
of~$G$.
\end{thm}
For $\K=\Q_p$, this is a classical result
by Wang \cite[Theorem~3.5]{Wan}.
In this case, $U_\alpha$ is always closed,
and $\alpha$ simply looks like a linear automorphism
in an exponential chart, which
facilitates to reduce the proof to
the case of linear automorphisms
of vector spaces. By contrast, automorphisms
of Lie groups over local fields
of positive characteristic
need not be linear in any chart
(cf.\ \cite{NOA}),
whence non-linear analysis cannot be~avoided.
\section{Examples of infinite-dimensional Lie groups}\label{seclie}
In this section, we discuss the
main examples of Lie groups over topological fields,
in parallel with some of the specific techniques of
non-linear functional analysis
needed to construct their Lie group structures.
In particular, we describe various
results concerning
non-linear mappings between function spaces.\vspace{-2mm}
\subsection*{Linear Lie groups}
\noindent
Among the easiest examples
of real or complex Lie groups are linear Lie groups,
i.e., unit groups
of unital Banach algebras (or other well-behaved
topological algebras)
and their Lie subgroups.
If $\K$ is a general topological field,
then a good class of topological algebras
to look at are the \emph{continuous inverse algebras},
i.e., unital associative topological $\K$-algebras~$A$
such that the group of units~$A^\times$ is open
and the inversion map $\eta \colon  A^\times \to A$, $a\mto a^{-1}$
is continuous. An elementary argument shows
(see \cite[Proposition~2.2]{ZOO}):
\begin{prop}
If $A$ is a continuous inverse algebra,
then the inversion map $\eta \colon A^\times \to A$
is smooth and thus $A^\times$ is a
Lie group.
\end{prop}
For example, $\K$ is a continuous inverse
algebra, and more generally
every finite-dimensional unital
associative $\K$-algebra when equipped
with the canonical vector topology ($\isom\K^d$),
see \cite[Proposition~2.6]{ZOO}.
If $A$ is a continuous inverse algebra
over~$\K$, then so is the matrix
algebra $M_n(A)$, for each $n\in \N$
(see \cite[Proposition~2.3]{ZOO}).
Furthermore, $A\tensor_\K \bL$ is a continuous
inverse algebra over~$\bL$,
for each finite extension $\bL$ of~$\K$,
by \cite[Corollary~2.8]{ZOO}.
If $A$ a continuous inverse algebra over
a locally compact topological field~$\K$,
and $K$ a compact $C^r$-manifold over~$\K$
(where $r\in \N_0\cup\{\infty\})$,
then also the algebra
$C^r(K,A)$ of $A$-valued
$C^r$-maps is a continuous inverse algebra,
with respect to pointwise operations
and the natural topology on this function space
described below
(see \cite[Proposition~5.7]{ZOO}).
For further examples over $\R$ or $\C$, see \cite{ALG}.\\[2.5mm]
Summing up, we always have a certain supply of continuous inverse
algebras over each topological field,
whose unit groups
provide a certain supply
of Lie groups.\vspace{-2mm}
\subsection*{Mapping groups and related constructions}
\noindent
The second widely studied class of infinite-dimensional
real Lie groups are the mapping groups, for example,
loop groups $C(\bS^1,G)$ and $C^\infty(\bS^1,G)$,
where $\bS^1$ is the unit circle and $G$ a finite-dimensional
real Lie group~(\cite{Mil}, \cite{PrS}).
The classical
constructions of mapping groups can be generalized
to a large extent to the case of Lie groups over topological fields.
For example, $C(K,G)$ can be made a Lie group
for each compact topological space~$K$
and Lie group $G$ over a topological
field~$\K$ (see \cite[Proposition~5.1]{ZOO}).
We now discuss
groups of differentiable
maps in more detail.
\cite[Proposition~5.1]{ZOO} subsumes:
\begin{prop}\label{mapgponcp}
Let $K$ be a compact
$($and hence finite-dimensional$)$
$C^k$-manifold over
a locally compact topological field~$\K$,
where $k\in \N_0\cup\{\infty\}$,
and $G$ a Lie group, modelled
on a topological $\K$-vector space~$E$.
Then there is a uniquely determined
smooth manifold structure on $C^k(K,G)$
making it a Lie group modelled on $C^k(K,E)$,
and such that
\[
\Phi\, \colon \; C^k(K,U)\to C^k(K,V), \quad
\gamma\mto \phi\circ \gamma
\]
defines a chart of $C^k(K,G)$ around $1$,
for some chart $\phi\colon G\supseteq U\to V\sub E$ of~$G$.
\end{prop}
Here
$C^k(K,U)=\{\gamma\!\in C^k(K,G)\colon \! \gamma(K)\sub U\}$
and
$C^k(K,V)=\{\gamma\!\in C^k(K,E)\colon \! \gamma(K)\sub V\}$,
which is an open subset of the topological
$\K$-vector space
$C^k(K,E)$ of $E$-valued $C^k$-maps
on~$K$.
The topology on the latter
is defined as follows.
\begin{defn}
Given topological $\K$-vector spaces $X$ and
$E$ over a topological field~$\K$,
and a $C^k$-map $\gamma \colon U\to E$
on an open subset $U\sub X$, where $k\in \N_0\cup\{\infty\}$,
we recursively define
$U^{[j]}:=(U^{[1]})^{[j-1]}$ and
$\gamma^{[j]} := (\gamma^{[1]})^{[j-1]}\colon
U^{[j]}\to E$
for each $j\in \N$ such that $j\leq k$
(with $U^{[0]}:=U$ and $f^{[0]}:=f$).
We equip $C^k(U,E)$
with the initial topology
with respect to the family
of mappings
\[
C^k(U,E)\to C(U^{[j]},E)
\,,\qquad \gamma\mto \gamma^{[j]}
\]
such that $j\in \N_0$ and $j\leq k$,
where $C(U^{[j]},E)$ carries the compact-open
topology.
If $M$ is a $C^k$-manifold
modelled on~$X$, with $C^k$-atlas $\cA$,
we equip $C^k(M,E)$ with the initial
topology with respect to
the family of maps
$C^k(M,E)\to C^k(V,E)$,
$\gamma\mto \gamma\circ \phi^{-1}$,
for $\phi\colon M\supseteq U\to V\sub X$ ranging through~$\cA$
(see \cite[\S4]{ZOO} for details).
\end{defn}
This is a very natural definition,
which provides spaces with the expected
completeness, metrizability
and local convexity properties
in relevant situations,
and which produces the conventional
topologies in the
real locally convex case~\cite[Proposition~4.19]{ZOO}.\\[2.5mm]
The proof of Proposition~\ref{mapgponcp}
is based on the following result
(and variants
for mappings between spaces of $C^k$-maps,
\cite[Proposition~4.20
and Corollary~4.21]{ZOO}).
It ensures smoothness
of the relevant non-linear mappings between
spaces of smooth functions:
\begin{prop}\label{onlysmth}
Let $K$ be a compact smooth manifold
over a locally compact topological field~$\K$,
$E$ and $F$ be topological $\K$-vector spaces,
and $U\sub E$ be open.
Then
\[
C^\infty(K,f)\, \colon \; C^\infty(K,U)\to C^\infty(K,F)\,,\qquad
\gamma\mto f\circ \gamma
\]
is a smooth map, for each smooth map $f\colon U\to F$.
More generally,
\[
f_*\, \colon \; C^\infty(K,U)\to C^\infty(K,F)\,,\qquad
\gamma\mto (x\mto f(x,\gamma(x)))
\]
is smooth, for each smooth map $f\colon K\times U\to F$.
\end{prop}
Now the $C^\infty$-version of Proposition~\ref{mapgponcp}
easily follows.
For example,
assuming $U=U^{-1}$,
the inversion
map of $C^\infty(K,G)$
is smooth on $C^\infty(K,U)$
(considered as a smooth manifold
with global chart~$\Phi$),
by the following argument:
The inversion map
of $G$ restricts to a smooth
map $i\colon U\to U$,
corresponding to the smooth map
$j:=\phi\circ i\circ \phi^{-1}\colon V\to V$
in the local chart~$\phi$.
The restriction
of the inversion map of $C^\infty(K,G)$
to $C^\infty(K,U)$
is $C^\infty(K,i)$.
This map is smooth,
since $\Phi\circ C^\infty(K,i)\circ\Phi^{-1}
=C^\infty(K,j)$ is smooth
by Proposition~\ref{onlysmth}.
\begin{rem}
We mention that
an analogue of Proposition~\ref{onlysmth}
is available for mappings
between spaces of
compactly supported vector-valued
$C^k$-functions
on a paracompact
finite-dimensional manifold~$M$ over
a locally compact topological field~$\K$
(see \cite[Proposition~8.22
and Corollary 8.23]{ZOO}).
For each Lie group~$G$
modelled on a topological $\K$-vector space~$E$,
this facilitates to turn the ``test function group''
\[
C^k_c(M,G)\;=\;\{\gamma\in C^k(M,G)\colon \mbox{$\gamma(x)=1$
for all $x$ outside some compact set}\}
\]
of compactly supported $G$-valued $C^k$-maps
into a Lie group
modelled on $C^k_c(M,E)$, equipped with a suitable
vector topology
(see \cite[Proposition~9.1]{ZOO}).
\end{rem}
As a tool for the discussion of diffeomorphism
groups, it is useful to know that the \emph{weak direct product}
\[
{\textstyle \prod_{i\in I}^*}\, G_i\, :=\,\big\{(x_i)_{i\in I}\in
{\textstyle \prod_{i\in I}G_i}\, \colon \;
\mbox{$x_i=1$ for all but finitely many $i$}\big\}
\]
of a family $(G_i)_{i\in I}$
of Lie groups over a valued field~$\K$
can always be turned into a Lie group,
modelled on the direct sum $E:=\bigoplus_{i\in I}E_i$
of the respective modelling spaces,
equipped with the box topology
(see \cite[Proposition~7.1]{ZOO}).
Thus,
sets of the form $\bigoplus_{i\in I}
U_i:=E\cap \prod_{i\in I}U_i$
are taken as a basis of open $0$-neighbourhoods for~$E$,
where each $U_i$ is an open $0$-neighbourhood in~$E_i$.
The following result
concerning typical
non-linear mappings between direct sums
(\cite[Proposition~6.9]{ZOO})
is used to construct the Lie group
structure on weak direct products.
\begin{prop}\label{mapsdirsums}
Let $(E_i)_{i\in I}$
and $(F_i)_{i\in I}$
be families of topological $\K$-vector spaces
indexed by a set~$I$. Let
$k\in \N_0\cup\{\infty\}$,
and
$f_i\colon U_i \to F_i$
be a $C^k$-map,
for $i\in I$,
defined on an open $0$-neighbourhood
$U_i\sub E_i$.
Then
$\bigoplus_{i\in I}\, f_i\, \colon \,
\bigoplus_{i\in I}\,U_i \to
\bigoplus_{i\in I} \,F_i\,$,
$(x_i)_{i\in I}\mto (f_i(x_i))_{i\in I}$
is a $C^k$-map on the open subset
$\bigoplus_{i\in I}U_i$ of
$\bigoplus_{i\in I}E_i$.
\end{prop}
For $\K\in \{\R,\C\}$, analogous results
can be obtained using locally convex
direct sums~\cite{MEA}.\vspace{-2mm}
\subsection*{Diffeomorphism groups}
\noindent
Let $M$ be a finite-dimensional, paracompact smooth manifold
over a local field~$\K$.
We explain some ideas used
in~\cite[\S13]{ZOO} to
define a Lie group structure
on $\Diff(M)$, the group of all
$C^\infty$-diffeomorphisms of~$M$:
\begin{prop}
$\Diff(M)$ is a Lie group
modelled on
the space $C^\infty_c(M,TM)$
of compactly supported smooth
vector fields on~$M$.
\end{prop}
To construct the Lie group structure, one exploits
that $M$ is a disjoint union $M=\bigcup_{i\in I}B_i$
of balls, i.e., open subsets $B_i\sub M$ diffeomorphic to~$\bO^d$,
where
$\bO$ is the maximal compact subring of~$\K$
and $d$ the dimension of the modelling space of~$M$.
The main step (explained more closely
below) is to make each
$\Diff(B_i)$ a Lie group.
Then the weak direct product
$\prod_{i\in I}^*\Diff(B_i)$
can be made a Lie group,
as described above.
Here $\prod_{i\in I}^*\Diff(B_i)$
can be identified with a subgroup
of~$\Diff(M)$ in an apparent way.
In a third step, one verifies that $\Diff(M)$
can be given a Lie group structure
making $\prod_{i\in I}^*\Diff(B_i)$ an open subgroup.\\[2.5mm]
The following ``exponential law''
(covered by \cite[Lemma~12.1
and Proposition 12.2]{ZOO})
is essential (cf.\ also
\cite[Proposition~12.6]{ZOO}
for a variant for metrizable manifolds).
\begin{prop}\label{prop48}
Let $\K$ be a topological field,
$M$ and $N$ be smooth manifolds
modelled on topological $\K$-vector spaces,
and $E$ be a topological
$\K$-vector space.
\begin{itemize}
\item[\rm (a)]
For every smooth mapping $f\colon M\times N\to E$,
also the mapping
$f^\vee  \colon M\to C^\infty(N,E)$,
$f^\vee(x):=f(x,\sbull)$ is smooth.
\item[\rm (b)]
If $\K$ is locally compact and $N$ is finite-dimensional,
then a map $g\colon M\to C^\infty(N,E)$
is smooth if and only if
$g^\wedge\colon M\times N\to E$,
$g^\wedge(x,y):=g(x)(y)$
is smooth.
\end{itemize}
Furthermore,
the map
$C^\infty(M\times N, E)\to C^\infty(M,C^\infty(N,E))$,
$f \mto f^\vee$ is an isomorphism
of topological vector spaces in the situation of {\rm (b)},
with inverse $g\mto g^\wedge$.
\end{prop}
We now explain the essential first step of the above
construction of
the Lie group structure on diffeomorphism
groups
in more detail,
namely the construction
of the Lie group structure
on $\Diff(M)$ for $M=\bO^d$ a ball.
Then $P:=\Diff(M)$ is an open subset
of the topological vector space $C^\infty(M,\K^d)$
(as above).\\[2.5mm]
\emph{Smoothness of inversion}.
The inclusion map
$i \colon P\to C^\infty(M,\K^d)$,
$\gamma\mto\gamma$
being smooth, the second half of the exponential law
(Proposition~\ref{prop48}\,(b))
ensures that also
\[
f \colon P\times M\to \K^d\,,\quad
f(\gamma,x):=i^\wedge(\gamma,x):=i(\gamma)(x)=\gamma(x)
\]
is smooth, using that~$M$ is finite-dimensional.
Then $f_\gamma:=f(\gamma,\sbull)=\gamma$
for each $\gamma\in P$,
whence $f_\gamma'(x)=\gamma'(x)\in \GL(\K^d)$
for each $x\in M$.
Applying the inverse function with parameters
(Theorem~\ref{invparultra}) with the diffeomorphism
$\gamma$ as the parameter, we see that
\[
g\colon P\times M\to M\,,\quad
g(\gamma,x):=(f_\gamma)^{-1}(x)=\gamma^{-1}(x)
\]
is smooth.
Now the first half of the exponential law
(Proposition~\ref{prop48}\,(a)) shows that
$\Diff(M)\to C^\infty(M,\K^d)$,
$\gamma\mto g^\vee(\gamma):=g(\gamma,\sbull)=\gamma^{-1}$
is smooth.
But this is the inversion map of the group $\Diff(M)$.\\[2.5mm]
\emph{Smoothness of composition}.
Since $\Diff(M)$ is open in
$C^\infty(M,\K^d)$, we only need to show
that
$\Gamma\colon C^\infty(M,\K^d\times \K^d)\isom
C^\infty(M,\K^d)\times C^\infty(M,\K^d) \to
C^\infty(M,\K^d)$, $\Gamma(\gamma,\eta):=\gamma\circ\eta$
is smooth.
Because the evaluation map
$\ev\colon C^\infty(M,\K^d)\times M\to\K^d$,
$\ev(\gamma,x):=\gamma(x)$
is smooth by \cite[Proposition~11.1]{ZOO},
we deduce from the formula
\[
\Gamma^\wedge((\gamma,\eta),x)=\gamma(\eta(x))
=\ev(\gamma,\ev(\eta,x))\quad\mbox{for $\, \gamma,\eta\in
C^\infty(M,\K^d)$, $x\in M$}
\]
that $\Gamma^\wedge$ is smooth and hence also $\Gamma$,
by Proposition~\ref{prop48}\,(b).
\begin{rem}
For differentiability properties
of the composition map between spaces
of $C^k$-functions,
see \cite[\S11]{ZOO}.
Diffeomorphism groups
can also be found in \cite{Lud}.\vspace{-2mm}
\end{rem}
\subsection*{Direct limit groups}
\noindent
Consider an ascending sequence
$G_1\sub G_2\sub \cdots$
of finite-dimensional
Lie groups over a locally compact topological field~$\K$,
such that each inclusion map
$G_n\to G_{n+1}$ is a smooth immersion.
Then $G:=\bigcup_{n\in \N}G_n$
can be given a Lie group structure
modelled on the direct limit topological $\K$-vector space
$\dl\;E_n$\vspace{-.8 mm}
of the respective finite-dimensional
modelling spaces,
which makes $G$ the direct limit of the given directed
sequence in the
category of $\K$-Lie groups
and smooth homomorphisms~\cite{FUN}
\,(cf.\ also \cite{NRW} and \cite{DIR}).
{\footnotesize
{\bf Helge Gl\"{o}ckner}, TU Darmstadt, FB Mathematik AG 5,
Schlossgartenstr.{}~7,
64289 Darmstadt, Germany;\\
E-Mail: {\tt gloeckner@mathematik.tu-darmstadt.de}}
\end{document}